\documentclass[letterpaper, 11pt,  reqno]{amsart}
\pdfoutput=1

\usepackage[margin=1.2in,marginparwidth=1.5cm, marginparsep=0.5cm]{geometry}

\usepackage{amsmath,amssymb,amsthm,amsxtra, esint}

\usepackage{color}
\usepackage[implicit=true]{hyperref}
\usepackage{amssymb}
\usepackage{mathrsfs}

\usepackage{xsavebox}

\usepackage{upgreek}

\usepackage{cases}

\allowdisplaybreaks[2]

\usepackage{tikz}

\definecolor{gr}{rgb}   {0.,   0.69,   0.23 }
\definecolor{bl}{rgb}   {0.,   0.5,   1. }
\definecolor{mg}{rgb}   {0.85,  0.,    0.85}
\definecolor{yl}{rgb}   {0.8,  0.7,   0.}
\definecolor{or}{rgb}  {0.7,0.2,0.2}

\newtheorem{theorem}{Theorem} [section]

\newtheorem{proposition}[theorem]{Proposition}
\newtheorem{remark}[theorem]{Remark}

\newtheorem{definition}[theorem]{Definition}

\newcommand{\noi}{\noindent}
\newcommand{\Z}{\mathbb{Z}}
\newcommand{\R}{\mathbb{R}}
\newcommand{\T}{\mathbb{T}}

\let\P= \undefined
\newcommand{\P}{\mathbf{P}}

\newcommand{\E}{\mathbb{E}}
\newcommand{\PP}{\mathbb{P}}

\newcommand{\al}{\alpha}
\newcommand{\be}{\beta}
\newcommand{\dl}{\delta}

\newcommand{\nb}{\nabla}

\newcommand{\Dl}{\Delta}
\newcommand{\eps}{\varepsilon}

\newcommand{\g}{\gamma}
\newcommand{\G}{\Gamma}
\newcommand{\ld}{\lambda}

\newcommand{\s}{\sigma}
\newcommand{\Si}{\Sigma}

\newcommand{\cj}{\overline}
\newcommand{\dx}{\partial_x}
\newcommand{\dt}{\partial_t}
\newcommand{\dd}{\partial}

\renewcommand{\o}{\omega}
\renewcommand{\O}{\Omega}

\newcommand{\les}{\lesssim}

\newcommand{\jb}[1]
{\langle #1 \rangle}

\newcommand{\ind}{\mathbf 1}

\newcommand{\M}{\mathcal{M}}

\newcommand{\N}{\mathbb{N}}
\newcommand{\NN}{\mathcal{N}}

\newtheorem*{ackno}{Acknowledgements}

\newcommand{\ze}{\zeta}
\newcommand{\upze}{\upzeta}

\newcommand{\too}{\longrightarrow}

\newcommand{\I}{\mathcal{I}}

\newcommand{\D}{\mathcal{D}}

\newcommand{\vp}{\varphi}
\newcommand{\bP}{\mathbb{P}}

\newcommand{\sF}{\mathscr{F}}

\newcommand{\cX}{\mathcal{X}}

\numberwithin{equation}{section}
\numberwithin{theorem}{section}

\DeclareMathOperator{\HS}{HS}

\makeatletter
\@namedef{subjclassname@2020}{%
	\textup{2020} Mathematics Subject Classification}
\makeatother

\begin{document}
\baselineskip = 14pt

\title[On probabilistic ill-posedness]
{On probabilistic ill-posedness}

\author[T.~Oh and N.~Tzvetkov]
{Tadahiro Oh and Nikolay Tzvetkov}
%

%

\address{
Tadahiro Oh, School of Mathematics\\
The University of Edinburgh\\
and The Maxwell Institute for the Mathematical Sciences\\
James Clerk Maxwell Building\\
The King's Buildings\\
Peter Guthrie Tait Road\\
Edinburgh\\ 
EH9 3FD\\
 United Kingdom, 
 and 
 School of Mathematics and Statistics, Beijing Institute of Technology,
Beijing 100081, China
}
\email{hiro.oh@ed.ac.uk}

\address{
Nikolay Tzvetkov,
Ecole Normale Sup{\'e}rieure de Lyon\\
UMPA\\
UMR CNRS-ENSL 5669\\
46, all{\'e}e d'Italie\\ 
69364-Lyon Cedex 07\\ 
France}

\email{nikolay.tzvetkov@ens-lyon.fr}



\begin{abstract}
In this note, 
we  introduce 
an enhanced notion of  probabilistic well-posedness
for dispersive PDEs with random initial data
by imposing stability at the origin
as the amplitude of randomization tends to $0$.
We then use this notion
to  re-interpret 
recent works on ``beyond variance blowup''
for  dispersive PDEs, by the authors with their collaborators (2025, 2026), 
as probabilistic ill-posedness results.
By drawing an analogy 
to the failure of $C^k$-smoothness
of a solution map in the deterministic setting, 
we interpret variance blowup results
as mild probabilistic ill-posedness.

\end{abstract}


\maketitle

\tableofcontents


\section{Introduction}

\subsection{Background}

In seminal works \cite{BO94, BO96}, 
Bourgain 
initiated the probabilistic well-posedness study
of nonlinear dispersive PDEs with random initial data
in the context of 
the  invariant Gibbs measure problem
for the nonlinear Schr\"odinger equation (NLS):
\begin{align}
i \dt u + \Dl u = |u|^{p-1} u.
\label{NLS1}
\end{align}

\noi
In particular, 
in \cite{BO96}, 
Bourgain studied the Cauchy problem
for the (renormalized) cubic NLS 
(with $p = 3$) on the two-dimensional torus $\T^2 = (\R/2\pi \Z)^2$
with the 
following Gaussian random  initial data (for $(d, \al)  = (2, 1)$):\footnote{By convention, we endow
$\T^d$ with the normalized Lebesgue measure $ dx_{_{\T^d}} =  (2\pi)^{-d}dx$
such that we do not need to carry factors involving $2\pi$.
With an abuse of notation, we simply use $dx$ to denote the normalized Lebesgue measure in the following.}
\begin{align}
u_0 = u_0^\o (\al) = \sum_{n\in \Z^d} \frac{g_n(\o) }{\jb{n}^\al} e_n ,
\label{gauss1}
\end{align}

\noi
where $\al \in \R$, 
$\jb{n} = \sqrt{1+|n|^2}$, 
$e_n(x) = e^{in\cdot x}$, 
and 
$\{g_n\}_{n\in \Z^d}$ is a sequence of independent 
complex-valued standard Gaussian random variables on some probability space 
$(\Omega, \sF, \bP)$.
When $\al = 1$, the random function $u_0$ in \eqref{gauss1}
corresponds to the so-called massive Gaussian free field on $\T^d$
whose law serves as the base Gaussian measure for the $\Phi^k_d$-measure on $\T^d$.
It is easy to see  that 
 $u_0$ in \eqref{gauss1}
belongs almost surely  to~$W^{s, p}(\T^d) \setminus W^{\al - \frac d2, p}(\T^d)$ for any 
$s < \al - \frac d2$ and  $1\leq p \leq \infty$.
Thus, for $(d, \al) = (2, 1)$, 
the random function $u_0$ in \eqref{gauss1}
is almost surely a distribution of negative regularity, 
 belonging 
 to $H^s(\T^2) \setminus L^2(\T^2)$
 for $s < 0$.
 Recall that the range $s < 0$
  corresponds to the so-called scaling supercritical regime, 
 where the (renormalized) cubic NLS on $\T^2$ is (deterministically) ill-posed~\cite{OhNI}.
In order to overcome this difficulty, 
Bourgain considered the following first order expansion:
\begin{align}
u = z + v, 
\label{B1}
\end{align}

\noi
where $z(t) = e^{it \Dl} u_0$ denotes the random linear solution, 
and  studied the equation satisfied by 
the remainder term $v$:
\begin{align}
\begin{cases}
i \dt v + \Dl v = \NN(v+z)\\
v|_{t = 0} = 0, 
\end{cases}
\label{NLS2}
\end{align}

\noi
where $\NN(u) = \big(|u|^2 - 2\int_{\T^2} |u|^2dx\big) u $
denotes the renormalized cubic nonlinearity.
See also~\cite{McK, DPD}, 
where analogous first order expansions were considered.
Then, by combining the 
harmonic analytic approach to study dispersive PDEs
(namely, the Fourier restriction norm method introduced in \cite{BO93})
with stochastic analysis
(constructing the renormalized powers 
of the random linear solution $z$), 
Bourgain proved that~\eqref{NLS2} is almost surely locally well-posed.
It is worthwhile to note that in analyzing the term $z^2 v$
(with a proper renormalization and a complex conjugate on 
one of the factors), 
he viewed the multiplication of $z^2$
as a random operator, say $M$,  and estimated
the operator norm of $M$
by studying the Hilbert-Schmidt norm of $MM^*$, 
which can be regarded as a precursor to the random tensor estimate
due to Deng, Nahmod, and Yue \cite{DNY3}.
See \cite[Section 1]{SLO}
for an intuitive explanation behind the random tensor estimate.
See also \cite{Bring2, Kaneshiro}
for an alternative proof of the random tensor estimate
based on the non-commutative Khintchine's inequality.

In \cite{BT1}, N.~Burq and the second author
studied a well-posedness issue for the cubic nonlinear wave equation (NLW)
on a three-dimensional compact Riemannian manifold
with randomized initial data 
that are given as randomization of 
(the Fourier coefficients of)
 a given pair of functions of 
the scaling supercritical regularity;
see \cite[p.\,451]{BT1}.
Proceeding with the first order expansion (as in \eqref{B1})
and exploiting the almost sure averaging property of the random
linear solution, 
they established almost sure local well-posedness
of the cubic NLW below the scaling critical regularity.
In \cite{BO96, BT1}, 
a  non-deterministic point of view has allowed them to go beyond the limit
of deterministic analysis 
(establishing almost sure local well-posedness below a scaling critical regularity,
where a given equation is (deterministically) ill-posed), 
and these works
\cite{BO96, BT1}
set a foundational basis   
of a (somewhat informal) notion of {\it probabilistic local well-posedness}, 
namely, existence and uniqueness
of a solution with random initial data
with stability under a {\it suitable} class of perturbations; 
 we will discuss the latter point in detail in the next subsection.
Over the last fifteen years, 
probabilistic well-posedness
of dispersive equations, broadly interpreted
with random initial data and\,/\,or (additive) stochastic forcing, 
has attracted much attention and has been studied intensively;
see, for example, 
\cite{BT1, CO, BT3, BOP2, GKO,  FOW, OTh2, 
 STz1,  Bring0, GKO2, GKOT, DNY2, STz2,  DNY3, OOT1, Bring1, 
OOT2, 
STzX, BDNY, BCST, Zine}.
See  also surveys \cite{BOP4, Tzv1, DNY5} on the subject.
We note that this development
also led to a 
 recent rapid development 
on 
pathwise well-posedness of stochastic dispersive PDEs with multiplicative noises;
see
\cite{CLO2, CGLLO2, CLOO, CLO3, COZ}.

In \cite{DNY2}, Deng, Nahmod, and Yue introduced the notion of 
probabilistic scaling critical regularity, 
which roughly
 corresponds to the regularity
above which 
the second Picard  iterate\footnote{More precisely, the second Picard iterate minus the linear solution.
In the following, we simply refer to such a term as the second Picard iterate.}
 (see, for example,  \eqref{exp2x})
is smoother than the random linear solution
(or the stochastic convolution in the stochastically forced case).
In fact, 
in order to make things a priori computable, 
a certain simplification 
was introduced in \cite{DNY2} in calculating
a probabilistic scaling critical regularity, 
by restricting one's attention to ``high $\times \cdots \times$ high $\to $ high'' interactions.
This notion of probabilistic scaling criticality
has provided useful heuristics
for threshold regularities, regarding probabilistic local well-posedness
of some dispersive equations.
For example,  in~\cite{DNY3}, 
the authors introduced the theory 
of random tensors
and showed
that NLS~\eqref{NLS1} (for $p \in 2\N+1$ with a renormalization on the nonlinearity)
is almost surely locally well-posed
with respect to the Gaussian random initial data of the form
\eqref{gauss1}
in the full probabilistically scaling subcritical regime.

As emphasized in \cite{DNY5}, 
 the probabilistic scaling heuristics 
  provides only a guiding principle
in determining a threshold regularity for probabilistic local well-posedness.  
In fact, 
there are recent results
\cite{Forlano, OO, Liu}
on {\it variance blowup} 
for dispersive PDEs
(see \eqref{div1} and Remark~\ref{REM:var1}), 
showing that, due to  divergence 
of the variance of a certain (renormalized) stochastic term, 
the standard probabilistic local well-posedness theory, 
based on the first order expansion 
\cite{McK, BO96, DPD}
(see \eqref{B1}) or its higher order variant
\cite{BOP3, OPTz, GKO2, OWZ},  
breaks down 
before reaching the threshold regularities
predicted by the probabilistic scaling heuristics;
see 
\cite[Remark~1.8]{OO}
and  \cite[discussion after Proposition 1.5]{Liu}.
Let us also mention a recent work \cite{BCST2} which establishes probabilistic well-posedness beyond the
probabilistic scaling regularity.

Following the work \cite{Hairer}
on the fractional KPZ equation, 
the authors with their collaborators
\cite{LLOT, LLLOT}
recently 
investigated 
  a possible 
extension
of probabilistic well-posedness theory
for dispersive PDEs beyond variance blowup
by considering a renormalization 
on the frequency-truncated random initial data
$\P_{N} u_0$
(where $\P_N$ denotes the frequency projection on 
frequencies $\{|n|\le N\}$)
via a vanishing multiplicative constant
$\dl_N \to 0$ as the frequency truncation parameter $N$ tends to $\infty$;
see, for example,  \eqref{BBM2}.

Our main goal in this note is
the following; 
by 
imposing continuity at the origin, 
we  introduce an enhanced  notion of probabilistic well-posedness
(see Definition \ref{DEF:1}).
We then use it to  re-interpret these ``beyond variance blowup''
results \cite{LLOT, LLLOT}
as {\it probabilistic ill-posedness} results.
In Subsection~\ref{SUBSEC:SPDE}, 
we also consider 
stochastic (dispersive\,/\,parabolic) PDEs, 
where 
``beyond variance blowup'' results
 such as \cite{Hairer, LLOT, LLLOT} imply
strong instability in the small noise limit.
In Section \ref{SEC:3}, 
by drawing an analogy 
to the failure of $C^k$-smoothness
of a solution map in the deterministic setting, 
we interpret variance blowup results
as mild probabilistic ill-posedness.

\begin{remark}\rm

Related to the probabilistic well-posedness study, 
there 
are results 
on 
pathological behavior of solutions with low regularity (random) initial data; see
\cite{STz0, OOT, OOPTz}
for further details.

\end{remark}

\subsection{Enhanced probabilistic local well-posedness}
\label{SUBSEC:1.2}

In this subsection, we introduce an enhanced notion of probabilistic local well-posedness
by paying particular attention to stability properties
of random solutions.
See Definition \ref{DEF:1}.

The standard notion of (deterministic) well-posedness
\`a la Hadamard \cite{Hada, Hada2}
says the following.
Consider the following dispersive PDE on $\T^d$:\footnote{For simplicity
of presentation, we focus on $\T^d$ in the following.}
\begin{align}
\dt u + Lu = \NN(u), 
\label{disp1}
\end{align}

\noi
where 
$L$ denotes a dispersive linear operator
and $\NN(u)$ denotes 
a nonlinearity.
We say that~\eqref{disp1} is locally well-posed
in a Banach space $B$
if, given any 
bounded set $E \subset B$, 
there exists $T = T(E) > 0$
such that, given any $u_0 \in \E$, 
there exists 
a unique solution $u \in C([0, T];B) $ to \eqref{disp1}
with $u|_{t = 0} = u_0$
on the time interval $[0, T]$. 
Moreover, the solution map, sending 
initial data $ u_0 \in E$ to 
a solution $u \in C([0, T];B) $ is continuous.

\begin{remark}\label{REM:ball}\rm

Many dispersive PDEs enjoy scaling symmetry, 
which induces the notion of scaling (sub)criticality.
The definition of local well-posedness
stated above
essentially corresponds to 
that in the scaling-subcritical regime, 
where the local existence time depends only on the size
$\|u_0\|_B$
of initial data $u_0$.
Such uniformity on local existence time
also appears in the random data Cauchy theory discussed below, 
where uniformity is in terms
of the frequency truncation parameter
$N$ for the frequency-truncated random initial data $\P_N u_0$
(see ``(ii)~stability under frequency truncation'')
or 
of $0 < \dl \le 1$
for
the scaled random initial data $\dl u_0$
(see ``(iii)~continuity at the origin'' below).
See also Definition \ref{DEF:1}(iii').

\end{remark}

\begin{remark}\label{REM:Hada}
\rm

In the following, we assume that when $u_0 = 0$, 
the trivial function $u \equiv 0$ is 
 the unique solution to~\eqref{disp1}, 
 which is  satisfied in the case
of a monomial nonlinearity we consider below.

In this setting, 
continuous dependence of the solution map
implies the following
{\it continuity at the origin}:
\begin{align}
\begin{split}
\hspace{5mm} & \text{Given a sequence  $\{u_{0, n}\}_{n \in \N}\subset B$, tending to $0$ as $n \to \infty$,
let $u_n \in C([0, T]; B)$}
\\
& \text{be the solution to \eqref{disp1}
with $u_n|_{t = 0} = u_{0, n}$. 
Then, as $n \to \infty$, $u_n$ converges to } 
\\
& \text{the trivial solution $u\equiv 0$ in $C([0, T]; B)$.
}
\end{split}
\label{Hada1}
\end{align}

\noi
The 
continuity at the origin is of particular importance
since 
examples for proving ill-posedness 
are often constructed to violate
the continuity at the origin, 
starting with 
Hadamard's  fundamental example
\cite[pp.\,33--34]{Hada3}
(in the context of the Laplace equation viewed as a Cauchy problem).
See 
also Remark \ref{REM:Hada2}.

\end{remark}

\begin{remark}\rm

We note that the uniqueness
of a solution $u$ in general holds only in a subspace $X
\subset C([0, T];B) $.
If the uniqueness holds in the entire class
$C([0, T];B)$, then we say that~\eqref{disp1} is unconditionally locally well-posed
and that the solution $u$ is unconditionally unique;
see
 \cite{KATO}.
See also \cite{GKONF, OW2, GOT, GOSW}
for a further discussion on unconditional uniqueness
of solutions to (stochastic) dispersive PDEs.
In the following, we 
 assume that solutions are conditionally unique
(in a suitable subclass of 
$C([0, T];B)$)
and do not discuss
an (unconditional) uniqueness issue further.

\end{remark}

Let us now turn to probabilistic local well-posedness.
Given $\al \in \R$,
fix $s < \al - \frac d2$
such that 
the random initial data $u_0$ in~\eqref{gauss1} belongs almost surely to 
$H^s(\T^d)$. 
In the following, we assume that 
\eqref{disp1} is (deterministically) ill-posed
in $H^s(\T^d)$, 
while it is (deterministically) locally well-posed
in $H^\s(\T^d)$ for some $\s > s$.
In particular, we assume that the classical solution map,  
continuously sending initial data 
$u_0$ in $H^\s(\T^d)$ 
to a solution 
 $u \in C([0, T]; H^\s(\T^d))$ to \eqref{disp1}, 
does not extend 
to a well-defined continuous map on $H^s(\T^d)$.

In this setting, 
the term ``probabilistic local well-posedness'' 
of \eqref{disp1}
in $H^s(\T^d)$  usually means
that there exists a set $\Si \subset \O$ with 
$\PP(\Si) = 1$ such that, for each $\o \in \Si$, 
the following holds:

\smallskip

\begin{itemize}
\item[\bf (i)] \underline{\bf existence \& uniqueness.}
There exist $T = T_\o > 0$ and a unique solution $u = u^\o
\in C([0, T]; H^s(\T^d))$ 
to \eqref{disp1}
with $u^\o|_{t = 0} = u_0^\o$
on the time interval $[0, T]$.
Here, both existence and uniqueness
of a solution are interpreted in a delicate sense, 
after introducing a suitable solution ansatz; see Step 2 below.

\medskip

\item[\bf (ii)] \underline{\bf stability under frequency truncation.}
Given $N \in \N$, 
let 
$ u_N$
be the  (smooth) random solution  to \eqref{disp1} with the frequency-truncated
random initial data
$ u_N|_{t = 0} = \P_{N} u_0$.
Then, 
for each $\o \in \Si$, 
we have
\begin{align}
 u_N^\o \too u^\o
\label{disp3}
\end{align}

\noi
 in $C([0, T];H^s(\T^d))$ as $N \to \infty$, 
 where  $u$ 
 denotes the solution
to \eqref{disp1} with the random initial data $u|_{t = 0}=  u_0$, 
constructed in Part (i), 
with the (random) local existence time
 $T = T_\o > 0$ 
 (possibly multiplied by a deterministic constant factor).

Similar convergence holds if we consider the smooth random solutions $u_\eps$
to \eqref{disp1} with the mollified random initial data
$u_\eps^\o|_{t = 0} = \rho_\eps * u_0^\o$, 
where $\rho$ is a mollification kernel and $\rho_\eps(x) = \eps^{-d} \rho(\eps^{-1}x)$, $\eps > 0$.
Namely, $u_\eps^\o$ converges 
in $C([0, T];H^s(\T^d))$ to the solution $u^\o$ 
to \eqref{disp1} with  $u^\o|_{t = 0}=  u_0^\o$
as $\eps \to 0$.\footnote{Here, 
the claimed convergence holds in general
for a fixed sequence $\{\eps_j\}_{j \in \N}$
of positive numbers, tending to $0$, 
and the set $\Si$ of full probability depends
on the choice of the sequence $\{\eps_j\}_{j \in \N}$.
For general $\eps \to 0$ (along the continuum), 
one may instead  prove convergence in probability. }

\end{itemize}

\smallskip

\noi
We note that the stability property
of random solutions
sensitively depends on how we approximate given random initial data.
In fact, given random initial data $u_0$, 
it is possible to introduce approximating  initial data such that 
the corresponding solutions almost surely diverge in the limit
(thus not converging to 
the solution $u$ 
to \eqref{disp1} with  $u|_{t = 0}=  u_0$).
See~\cite{STz0, OOT} for a further discussion.


A standard approach to proving such probabilistic well-posedness
for \eqref{disp1}, 
based on the first order expansion 
 or its higher order variant, 
is to decompose
the 
classically ill-defined solution map
into two steps.

\smallskip

\begin{itemize}
\item[\bf Step 1:]
Given the random initial data $u_0$ in \eqref{gauss1}, 
we 
use stochastic analysis 
and
construct an enhanced data set $\Xi = \Xi(u_0)$:
\begin{align}
\Xi = \Xi (u_0)= \big(\Xi_1(u_0), \dots, \Xi_K(u_0)\big) 
\label{enh1}
\end{align}

\noi
consisting of a finite number of random distributions
(such as renormalized powers of the random linear solution, possibly 
under the Duhamel integral operator)
and random operators, 
and the space 
$ \cX_K$ for enhanced data sets
such that 
$\Xi(u_0) \in \cX_K$, almost surely.

For the sake of discussion, we assume that 
$\Xi_1$ is given by the random linear solution:
\begin{align}
\Xi_1(t) = e^{-tL} u_0
\label{disp1a}
\end{align}

\noi
in the following.

\medskip

\item[\bf Step 2:]
Given a (deterministic) enhanced data set $\Xi = (\Xi_1, \dots, \Xi_K) \in \cX_K$, 
write $u$ as
\begin{align}
u = \Xi_1 + \dots + \Xi_k + v
\label{exp1a}
\end{align}
for some $k \in \{1, \dots, K\}$ 
and consider the equation satisfied by the remainder term~$v$:\footnote{with a possible renormalization on the nonlinearity. In this discussion, we ignore the issue of renormalization.}
\begin{align}
\begin{cases}
\dt v + Lv = \NN(v + \Xi_1 + \dots + \Xi_k ) 
- (\dt + L) (\Xi_2 + \dots + \Xi_k ) \\
v|_{t = 0} = 0.
\end{cases}
\label{disp2}
\end{align}

\noi
By deterministic analysis (often via a contraction argument), 
we construct a continuous map:
\begin{align}
\Xi \in \cX_K
\longmapsto v \in C([0, T]; H^{\s}(\T^d))
\label{disp4}
\end{align}

\noi
for suitable $\s > s$
(recall that we assumed  that the original equation \eqref{disp1}
is (deterministically) locally well-posed in $H^{\s}(\T^d)$).
Then, we construct a solution $u$ to \eqref{disp1}
via the relation~\eqref{exp1a}.

\end{itemize}

\smallskip

\noi
For example, 
under 
the first order expansion 
$u = \Xi_1 + v = e^{-tL}u_0  + v$ (as in \eqref{B1}),  
the remainder term $v = u - \Xi_1$ satisfies the following equation:
\begin{align*}
\dt v + Lv = \NN(v+\Xi_1) 
 \quad \text{with \ $v|_{t = 0} = 0$}.
\end{align*}

\noi
Let 
 \begin{align}
 \Xi_2 =  \I\big(\NN(\Xi_1) \big)
\label{exp2x} 
 \end{align}
denote the second Picard iterate,
where $\I$ denotes the Duhamel integral operator:
\begin{align}
\I(F)(t) = \int_0^t e^{-(t-t')L} F(t') dt'.
\label{disp4b}
\end{align}

\noi
Then, the second order expansion:
\begin{align}
u = \Xi_1 + \Xi_2 +v = \Xi_1 + \I\big(\NN(\Xi_1)\big) + v
\notag 
\end{align}

\noi
yields that 
the remainder term $v = u - \Xi_1 - \Xi_2$ satisfies the following equation:
\begin{align*}
\dt v + Lv = \NN(v+\Xi_1 + \Xi_2) - \NN(\Xi_1)
 \quad \text{with \ $v|_{t = 0} = 0$}.
\end{align*}


\begin{remark}\label{REM:sto1}\rm

In Step 1,  we construct the enhanced data set $\Xi(u_0)$
in \eqref{enh1}
by first constructing
the enhanced data set  $\Xi(\P_N u_0)$
for  the frequency-truncated initial data 
$\P_N u_0$
and then 
by establishing 
the following 
almost sure convergence:
\begin{align}
\Xi(\P_N u_0)\,  \too \,  \Xi(u_0)
\label{disp5}
\end{align}

\noi
in $\cX_K$, 
as $N \to \infty$.
Then, 
the almost sure convergence \eqref{disp3}
of $u_N$ to $u$ 
follows from 
the continuity of the map 
\eqref{disp4} in Step 2
with \eqref{disp5}.
Similar convergence holds
if we replace the frequency truncation
by mollification.

\end{remark}

For  simplicity of presentation, let us focus on 
 the case of a monomial nonlinearity: 
 \begin{align}
  \NN(u) = u^p\quad \text{(possibly with derivatives)}
\label{mono1}  
 \end{align}
 
 \noi
for an integer  $p \ge 2$ 
 in~\eqref{disp1}
 (and $u^{p_1} \cj u^{p_2}$ with $p = p_1 + p_2$
in the complex-valued setting (possibly with derivatives)).
In this case, 
we typically have the following;\footnote{Even 
if $\Xi_k(\dl u_0) $ is not homogeneous in $\dl$, 
we still have 
$\Xi_k(\dl u_0)  =  \dl^{\ld_k} O(\Xi_k (u_0))$
in general when measured in some space-time norm 
or operator norm, from which 
 the convergence \eqref{disp6} follows.}
for each  $k = 1, \dots, K$, 
there exists 
$\ld_k \in \N$ such that 
\begin{align}
\Xi_k(\dl u_0) = \dl^{\ld_k} \Xi_k (u_0) 
\label{disp5a}
\end{align}

\noi
for any $\dl > 0$.
Namely, 
the enhanced data set 
 $\Xi(\dl  u_0)$
associated with  
the scaled random initial  data 
$\dl u_0$
converges almost surely to~$0$
in $\cX_K$
as $\dl \to 0$:\footnote{In the current monomial case, 
we impose that $\Xi (0) = (0, \dots, 0)$.}
\begin{align}
\Xi(\dl u_0)= \big(\Xi_1(\dl u_0), \dots, \Xi_K(\dl u_0)\big)\,  \too \, \Xi (0) = (0, \dots, 0).
\label{disp6}
\end{align}

\noi
This discussion leads to  the following 
continuity
 at the origin.

\smallskip

\begin{itemize}
\item[\bf (iii)] \underline{\bf continuity at the origin.}
Given $\dl > 0$, 
let 
$ u^\dl$
be the   random solution  to~\eqref{disp1} with the 
scaled random  initial data
$ u^\dl|_{t = 0} = \dl u_0$, 
given by the expansion (see \eqref{exp1a}):
\begin{align}
u^\dl  = \Xi_1 (\dl u_0) + \dots + \Xi_k(\dl u_0) + v^\dl, 
\notag 
\end{align}

\noi
where $v^\dl$ is an almost sure solution 
to \eqref{disp2} with $\Xi = \Xi(\dl u_0)$.
Note that, in view of~\eqref{disp5a} and~\eqref{disp6}, 
we can choose the random local existence  time $T = T_\o$
(which is positive for each $\o \in \Si$, 
where $\Si$ is as in (i) and (ii))
to be independent of $0  < \dl \le 1$.

Combining \eqref{disp6} with Step 2 above
(in particular, the continuity of the map in~\eqref{disp4}), 
we have, for each $\o \in \Si$, 
\begin{align}
 u^{\dl, \o} \too  0
\label{disp7}
\end{align}

\noi
 in $C([0, T];H^s(\T^d))$ as $\dl \to 0$, 
 where the trivial limit $0$ denotes the function $u \equiv 0$
 which is a solution to \eqref{disp1} with the zero initial data.

\end{itemize}

\smallskip

\noi
See 
\cite{OQ} for a  discussion
on the large deviation principle (in a mild form)
related to the convergence \eqref{disp7}.

\begin{remark}\rm 
In   the paracontrolled approach
\cite{GKO2, OOT1, Bring1, OOT2},
we further split the remainder term $v$ in \eqref{exp1a}
into the singular part and the smoother part.
We note that even for such a paracontrolled approach, 
the properties (i), (ii), and  (iii)
still hold.
We also point out that 
the properties (i), (ii), and  (iii)
also hold
for solutions constructed via
more involved 
approaches
such as those 
 in \cite{Oh1, Oh2, OTzW, Bring0}, 
the
random averaging operators \cite{DNY2}, 
and the theory of random tensors
\cite{DNY3}.
A similar comment applies to the property (iii') 
stated in Definition~\ref{DEF:1} below.

For 
non-monomial nonlinearities considered in 
\cite{STz00, ORSW1, ORSW2, ORW, GHOZ, Zine, FZ}, 
a separate discussion is needed, but we do not discuss the matter further 
in this paper.

\end{remark}

\begin{remark}\rm 
Unlike the deterministic local well-posedness theory \`a la Hadamard
described above, 
continuous dependence on initial data does not hold
in general 
for the random data Cauchy theory, 
since the random data (for example, $u_0 = u_0(\al)$ in \eqref{gauss1})
is defined only almost everywhere 
 $H^s(\T^d)$ in a suitable sense, {\it not} everywhere in $H^s(\T^d)$.
Nonetheless, in~\cite{BT3}, 
N.\,Burq and the second author introduced the notion
of {\it probabilistic continuous dependence}, 
utilizing the conditional probability 
 that two random solutions
with the random initial data 
distributed by the same law (say, by \eqref{gauss1})
being apart by the distance $\eps > 0$, 
given that two initial data are within the distance $\dl$, 
tends to $0$ 
as $\dl \to 0$
for each {\it fixed} $\eps > 0$; 
 see \cite[Theorem 3]{BT3}.
See also \cite{Poc}.

The continuity at the origin (iii) stated above
is of different nature than the probabilistic continuous dependence introduced in \cite{BT3}
in the sense that it is a pathwise statement (namely, for each $\o \in \Si$ with $\PP(\Si) = 1$)
without relying on the conditional probability as in \cite{BT3}.
A similar comment applies to Definition \ref{DEF:1}\,(iii').
\end{remark}

Before proceeding further, we note the following.
Given   any  sequence $\{\dl_N\}_{N \in \N}$
of positive numbers, tending to $0$  as $N \to \infty$, 
it follows from writing 
$ \Xi_k (\dl_N \P_N u_0)
= \dl_N^{\ld_k} \Xi_k (\P_N u_0)$
as
\begin{align}
 \Xi_k (\dl_N \P_N u_0)
= \dl_N^{\ld_k} \big(\Xi_k ( \P_N u_0)-
\Xi_k ( u_0)\big)
+ 
\dl_N^{\ld_k} 
\Xi_k ( u_0)
\label{disp7a}
 \end{align}

\noi
with  \eqref{disp5} 
that 
\begin{align}
\Xi(\dl_N \P_N u_0)\,  \too \,  \Xi (0) = (0, \dots, 0)
\label{disp9}
\end{align}

\noi
in $\cX_K$, 
almost surely (namely for each  $\o \in \Si$), 
as $N \to \infty$.
This motivates the following enhanced notion of probabilistic local well-posedness.

\begin{definition}\label{DEF:1}\rm

Given $\al \in \R$,
fix $s < \al - \frac d2$ as above.
We say that 
probabilistic local well-posedness holds
for 
the nonlinear dispersive PDE \eqref{disp1}
with respect to the random initial data $u_0$ of the form \eqref{gauss1}, 
if
 there exists 
a set $\Si \subset \O$ with $\PP(\Si) = 1$ such that, for each $\o \in \Si$, 
the following properties hold:

\smallskip

\begin{itemize}
\item[\bf (i)] {\bf existence \& uniqueness},

\smallskip
\item[\bf (ii)]
{\bf  stability under the frequency truncation}

\end{itemize}

\smallskip

\noi
as stated above, and 

\smallskip

\begin{itemize}
\item[\bf (iii')]
\underline{\bf continuity at the origin.}
Let  $\{\dl_N\}_{N \in \N}$
be 
 a sequence  of positive numbers tending to $0$
 as $N \to \infty$.
Then, denoting by  
$ u_N$
 the  (smooth) random solution  to \eqref{disp1} with the scaled frequency-truncated
 random initial data
$ u_N|_{t = 0} = \dl_N \P_{N} u_0$, 
we have, for each $\o \in \Si$, 
\begin{align}
 u_N^\o \too 0 
\label{disp10}
\end{align}

\noi
 in $C([0, T];H^s(\T^d))$ as $N \to \infty$, 
 where  
 $T = T_\o > 0$ denotes the (random) local existence time
 from Part (i) (possibly multiplied by a deterministic constant factor).

\end{itemize}


\end{definition}

In view of the continuity of the map in \eqref{disp4}
and the convergence 
\eqref{disp9}
of the enhanced data set 
$\Xi(\dl_N \P_N u_0)$, 
we see that standard probabilistic local well-posedness
results in the literature satisfy
the continuity at the origin~(iii').
As such, 
the condition (iii') in 
Definition~\ref{DEF:1} is a natural one to impose
as part of criteria for probabilistic local well-posedness.

In the next section, we present several examples of probabilistic ill-posedness
results with respect to 
this enhanced notion of probabilistic local well-posedness
introduced 
in 
Definition~\ref{DEF:1}.

\begin{remark}\label{REM:Hada2}\rm
As pointed out in Remark \ref{REM:Hada}, 
in proving ill-posedness in the deterministic setting, 
one often aims to violate 
the continuity at the origin \eqref{Hada1}.
In an analogous manner, 
we 
establish  probabilistic ill-posedness
by violating 
the continuity at the origin 
\eqref{disp10} in
Definition~\ref{DEF:1}\,(iii').
See Section \ref{SEC:2}.
\end{remark}

\begin{remark}\label{REM:CM} \rm

Fix a deterministic function $v_0 \in H^\s(\T^d)$.
(Recall that, by our assumption, \eqref{disp1}
is locally well-posed in $H^\s(\T^d)$.) 
Given a sequence   $\{\dl_N\}_{N \in \N}$
 of positive numbers tending to $0$
 as $N \to \infty$, 
 let $u_N$ be the solution to \eqref{disp1}
with initial data $ u_N|_{t = 0} = v_0 + \dl_N \P_N  u_0$.
Then, with the expansion \eqref{exp1a}, 
the remainder term $v_N = u_N - 
\big(\Xi_1(\dl_N \P_N  u_0) + \dots + \Xi_k(\dl_N \P_N  u_0)\big)$
satisfies 
the following equation:
\begin{align}
\begin{cases}
\dt v_N + Lv_N = \NN(v_N + \Xi_{1, N} + \dots + \Xi_{k, N} ) 
- (\dt + L) (\Xi_{2, N} + \dots + \Xi_{k, N} ) \\
v_N|_{t = 0} = v_0, 
\end{cases}
\notag
\end{align}

\noi
where $\Xi_{j, N} = 
\Xi_j(\dl_N \P_N  u_0)$
with 
$\Xi_{1, N} (t)= 
\Xi_1(\dl_N \P_N  u_0)(t) = \dl_N e^{-tL} \P_N u_0$
as in 
\eqref{disp1a}.
Then, together with 
\eqref{disp9}, 
a slight modification of the almost sure local well-posedness
argument described above 
yields
the following continuity 
at a general deterministic function $v_0 \in H^\s(\T^d)$.

\smallskip

\begin{itemize}
\item[\bf (iii'')]
\underline{\bf continuity at a general deterministic function $v_0 \in H^\s(\T^d)$.}
Given a sequence   $\{\dl_N\}_{N \in \N}$
 of positive numbers tending to $0$
 as $N \to \infty$, 
 let $u_N$ be the solution to \eqref{disp1}
with initial data $ u_N|_{t = 0} = v_0 + \dl_N \P_N  u_0$.
Then,  for each $\o \in \Si$, we have 
\begin{align*}
 u_N^\o \too  u
\end{align*}

\noi
 in $C([0, T];H^s(\T^d))$ as $N \to \infty$, 
 where $u$ is the solution to \eqref{disp1}
 with the deterministic initial data $u|_{t= 0} = v_0 \in H^\s(\T^d)$.

\end{itemize}

\end{remark}

\section{``Beyond variance blowup'' as probabilistic ill-posedness}
\label{SEC:2}

In this section, we go over several examples
of recent ``beyond variance blowup'' results
and interpret them as  probabilistic ill-posedness results.
While it is immediate from Definition~\ref{DEF:1} and Egoroff's theorem, 
the following proposition 
provides an important criterion for probabilistic ill-posedness.

\begin{proposition}\label{PROP:ill1}

Given $\al \in \R$, let $u_0 = u_0(\al)$
be the Gaussian random function in \eqref{gauss1}.
Suppose that there exists a sequence $\{\dl_N\}_{N \in \N}$
of positive numbers, tending to $0$ as $N \to \infty$, 
such that 
\begin{align}
\limsup_{N \to \infty} \PP\Big( \|u_N\|_{C_{T} H^s_x} > \ld\Big)  >  0 
\label{disp11}
\end{align}

\noi
for some $\ld > 0$.
Here, 
$ u_N$
denotes 
 the  (smooth) random solution  to 
 the nonlinear dispersive PDE 
 \eqref{disp1} with the scaled frequency-truncated
 random initial data
$ u_N|_{t = 0} = \dl_N \P_{N} u_0$
on the time interval $[0, T]$, 
where $T = T_\o > 0$ denotes the random local existence time.
Then, 
the probabilistic Cauchy problem \eqref{disp1} 
with the Gaussian random initial data
$u_0 = u_0(\al)$ in~\eqref{gauss1}
is  ill-posed
for this particular value of $\al$.

\end{proposition}

The condition  \eqref{disp11}
violates 
``the continuity at the origin''  in Definition \ref{DEF:1}\,(iii'), 
resulting in probabilistic ill-posedness of \eqref{disp1}.
The ``beyond variance blowup'' results mentioned
in the subsequent subsections
in fact state that $u_N$ converges
in law 
to a non-trivial solution, verifying the condition~\eqref{disp11};
see \eqref{BBM5} and \eqref{NLW4}.
We, however,  point out that 
probabilistic ill-posedness
of \eqref{disp1}
follows as long as  $u_N$ (up to a subsequence)
{\it remains bounded away from~$0$} in a suitable sense, 
and thus no convergence
of $u_N$ (to a non-trivial limit) is required.

In the next two subsections, we go over several examples
of probabilistic ill-posedness.

\subsection{Benjamin-Bona-Mahony equation}
\label{SUBSEC:2.1} 

We first consider the following Benjamin-Bona-Mahony equation (BBM)
 on  the circle $\T$:
\begin{align}
\begin{cases}
\dt u-\partial_{txx}u+\partial_{x}u
+\partial_{x}(u^{2})  = 0  \\
u|_{t = 0} = u_0, 
\end{cases}
\quad  (t,x) \in \R\times \T,
\label{BBM0}
\end{align}

\noi
where $u$ is a real-valued unknown.
With 
$D=-i\partial_{x}$, define the operator 
$\vp(D)$ by 
\begin{align}
\vp(D)=
- (1-\partial_{x}^{2})^{-1}\partial_{x}.
\notag 
 \end{align}


\noi
Then, we can rewrite \eqref{BBM0}
 as
\begin{align}
\begin{split}
\dt u 
& = \varphi(D)  u+ \NN(u) \\
:\! & = \varphi(D)  u+\varphi(D) (u^2). 
\end{split}
\label{BBM1}
\end{align}

\noi
In the following, we  use~\eqref{BBM0} and \eqref{BBM1}
interchangeably.

The BBM equation \eqref{BBM0}
was 
 proposed  in \cite{BBM} as an alternative model to the Korteweg-de Vries  equation
 (see \eqref{KDV1} below),
 and it is known to be globally well-posed
in $H^s(\T)$ for $s \ge 0$, 
while it is ill-posed for $s < 0$;
see \cite{BT09, RO10,  PA11,  BD16, Forlano}.
In an effort to extend the well-posedness theory to Sobolev spaces
of negative regularity, Forlano \cite{Forlano}
studied the well-posedness issue of~\eqref{BBM1} with  the 
Gaussian random initial data $u_0 = u_0(\al)$ in \eqref{gauss1}
(conditioned that  $g_{-n} = \cj{g_n}$ in the current real-valued setting).
We note that 
for $\al > \frac 12$, the Gaussian random initial data $u_0(\al)$
in \eqref{gauss1} almost surely belongs to $L^2(\T)$, 
to which deterministic global well-posedness of \eqref{BBM0} applies, 
and thus we restrict our attention to $\al \le \frac 12$ in the following.
By considering the first order expansion 
$u = \Xi_1(u_0) + v$, 
where $\Xi_1(u_0)(t) = e^{t\varphi(D)}u_0 $ denotes the random linear solution, 
and studying the equation satisfied by the remainder term
$v$, 
he 
showed that, when $\frac 14 < \al \le \frac 12$, 
BBM~\eqref{BBM1} 
(hence~\eqref{BBM0})
is almost surely locally well-posed
with respect to the Gaussian random initial data $u_0$ in \eqref{gauss1}.

In the same paper \cite{Forlano}, Forlano exhibited the 
following  variance blowup for $\al \le \frac 14$.
Let 
$\Xi_2(f)$ denote the second Picard iterate defined by 
 \begin{align}
 \Xi_2 (f)(t)=  
 \int_0^t 
e^{ (t- t') \varphi(D)} \NN(e^{t'\varphi(D)} f ) dt'.
\label{Pi2}
 \end{align}

\noi
Then, for $\al \le \frac 14$, 
we have 
\begin{align}
\lim_{N \to \infty} 
\E\Big[|\jb{  \Xi_2(\P_N u_0)(t), \psi}|^2\Big] = \infty
\label{div1}
\end{align}

\noi
for any $t > 0$
and any non-zero test function $\psi \in \D(\T) = C^\infty(\T)$, 
where 
$\jb{f, g}$ 
denotes the $\D'(\T)$-$\D(\T)$ duality pairing.

\begin{remark}\label{REM:var1}\rm

As in \cite{Hairer, LLOT}, we use the term  ``variance blowup''
to describe the situation, 
where (i)~the variance of a certain stochastic term diverges (see \eqref{div1}
for the BBM case) at a threshold
{\it but} (ii)~the analytical framework 
for proving almost sure local well-posedness
still continues to hold beyond the threshold.
See  \cite[Remark~1.2]{LLOT}
for a further discussion.

\end{remark}

\begin{remark}\rm
We recall that 
BBM \eqref{BBM1} with the Gaussian random initial data $u_0 = u_0(\al)$
in \eqref{gauss1} is 
critical with respect to the 
probabilistic scaling 
when 
$\al_\text{crit} = 0$.
On the other hand, 
the variance blowup \eqref{div1} occurs at $\al = \frac 14$
before reaching 
the probabilistic scaling critical threshold
$\al_\text{crit} = 0$.
This is due to the fact that the variance blowup \eqref{div1}
is based on the 
``high $\times$ high $\to$ low'' interaction
(and not 
on the ``high $\times$ high $\to$ high'' interaction, 
from which 
the probabilistic scaling critical threshold
$\al_\text{crit} = 0$
is determined).

\end{remark}

In  \cite{LLOT}, 
 the authors with G.\,Li  and J.\,Li
studied BBM 
beyond the  variance blowup, 
following Hairer's  work \cite{Hairer}
on the  fractional KPZ equation.
Fix $ \al \le \frac 14$.
Given $N \in \N$, we introduce 
the multiplicative renormalization constant 
  $\dl_{\al, N}$ by setting 
\begin{align}
\dl_{\al, N}= 
\bigg( \sum_{|n|\le N }
\frac{2}{\jb{n}^{4\al}}\bigg)^{-\frac 14}
\sim 
\begin{cases}
(\log \jb{N})^{-\frac 14} & \text{for } \al = \frac 14,\\
\jb{N}^{\al-\frac 14}
  &  \text{for }   \al<\frac 14.
\end{cases}
\label{CN}
\end{align}

\noi
Note that 
\begin{align}
\dl_{\al, N} \too 0\quad \text{as}\quad N \too \infty.
\notag 
\end{align}

\noi
Then,   consider  BBM 
\eqref{BBM0}
with the renormalized (frequency-truncated) random initial data:
\begin{align}
\begin{cases}
\dt u_N-\partial_{txx}u_N+\partial_{x}u_N
+\partial_{x}(u_N^{2})  = 0  \\
u_N|_{t=0} = \dl_{\al, N} \P_{N} u_0.
\end{cases}
\label{BBM2}
\end{align}

\noi
In \cite{LLOT}, 
the authors established the following ``beyond variance blowup'' result;
as $N \to \infty$, 
the (smooth) solution $u_N$
to \eqref{BBM2}
converges in law
to the solution $u$ to 
the following stochastic BBM  forced by 
the derivative of
a spatial white noise
 $\ze$:
\begin{align}
\begin{cases}
\dt u-\partial_{txx}u+\partial_{x}u
+\partial_{x}(u^{2})  = \dx \ze  \\
u|_{t = 0} = 0
\end{cases}
\label{BBM4d}
\end{align}

\noi
in $C(\R_+; H^{- \frac 14 - \eps}(\T))$ for any $\eps > 0$.
We point out that this is a ``central limit theorem''-type result.

The main observation 
for this convergence result
is the following; 
due to our choice~\eqref{CN} of the renormalization constant $\dl_{\al, N}$, 
the quadratic term in \eqref{BBM2} converges
to the linear term 
in~\eqref{BBM4d}.
More precisely, 
as $N \to \infty$, the second Picard iterate 
$\Xi_2(\dl_{\al, N} \P_N u_0)$ for \eqref{BBM2}, 
where  $\Xi_2 $ is as in \eqref{Pi2}, 
converges in law to 
the linear solution to~\eqref{BBM4d}: 
\begin{align*}
Z(t)=  
- 
 \int_0^t e^{ (t- t') \varphi(D)} \varphi(D)\zeta  dt'
 \end{align*}

\noi
in $C([0, T]; W^{\frac 12 - \eps, \infty}(\T))$
for any $T>0$ and $\eps > 0$, 
where 
$\ze$ denotes the spatial white noise as in \eqref{BBM4d}.
Once this convergence-in-law 
of $\Xi_2(\dl_{\al, N} \P_N u_0)$ to $Z$
is established, 
by applying the Skorokhod representation theorem \cite[Chapter 31]{Bass}
and working with the second order expansion for \eqref{BBM2}, 
we can prove
 almost sure convergence of 
$u_N$ to $u$ on a new probability space, 
which in turn implies convergence in law of $u_N$ to $u$
on the original probability space;
see \cite{LLOT} for the details.
Then, by the portmanteau theorem, 
there exists small $\ld> 0$ such that 
\begin{align}
\begin{split}
\limsup_{N \to \infty} 
\PP\Big(\| u_N\|_{C_T H^{-\frac 14- \eps}_x} > \ld  \Big)
& \ge 
\liminf_{N \to \infty} 
\PP\Big(\| u_N\|_{C_T H^{-\frac 14- \eps}_x} > \ld  \Big)\\
& \ge \PP\Big(\| u\|_{C_T H^{-\frac 14- \eps}_x} >\ld\Big)
> 0 
\end{split}
\label{BBM5}
\end{align}

\noi
for any $T > 0$, 
where 
the last inequality follows from 
the fact that the solution $u$ to~\eqref{BBM4d}
is not trivial (namely,  $u \not \equiv 0$), almost surely.
Hence, from Proposition \ref{PROP:ill1}, 
we conclude the following probabilistic ill-posedness
result for BBM \eqref{BBM0} on $\T$.

\begin{theorem}
\label{THM:1}

Let $\al\leq \frac14$.
Then, 
the probabilistic Cauchy problem 
for BBM \eqref{BBM0}  on $\T$
with the Gaussian random initial data
$u_0 = u_0(\al)$ in \eqref{gauss1}
is  ill-posed
in the sense of Definition~\ref{DEF:1}.

\end{theorem}

\subsection{Quadratic nonlinear wave equation}
\label{SUBSEC:2.2}

Consider the quadratic nonlinear wave equation (NLW)
on $\T^2$:
\begin{align}
\begin{cases}
\dt^2 u + (1- \Dl) u = \NN(u)\\ 
(u, \dt u)|_{t = 0} = (u_0, u_1), 
\end{cases}
\quad  (t,x) \in \R\times \T^2,
\label{NLW1}
\end{align}

\noi
where  
$\NN(u) = \, :\!u^2\!:$ denotes the Wick renormalization of
the quadratic nonlinearity  $u^2$ and 
 the random initial data 
$(u_0, u_1) = (u_0^\o(\be), u_1^\o(\be)) $
is given 
by\footnote{Here,  
we make a different parameter choice than \eqref{gauss1}
(namely, $\be = 1-\al$)
to match the notation in~\cite{OO, LLLOT}.} 
\begin{align}
(u_0, u_1)   =
\bigg( \sum_{n\in \Z^2} \frac{g_n(\o) }{\jb{n}^{1-\be}} e_n,
\sum_{n\in \Z^2} \jb{n}^\be h_n(\o)  e_n \bigg).
\label{gauss2}
\end{align}

\noi
Here, 
$\{g_n, h_n \}_{n\in \Z^2}$ is a family  of independent 
complex-valued standard Gaussian random variables
conditioned that $g_{-n} = \overline{g_n}$
and $h_{-n} = \overline{h_n}$, $n \in \Z^2$.
Then, a straightforward modification of the 
argument in \cite[Proposition 1.6]{OO}
(carried out in the context of the stochastic NLW 
\eqref{SNLW1})
shows that 
variance blowup on the second Picard iterate
occurs
for  $\be \ge \frac 12$:
\begin{align}
\lim_{N \to \infty} 
\E\Big[\big|\big\langle  \Xi_2\big((\P_N u_0, \P_N u_1)\big)(t), \psi \big\rangle\big|^2\Big] = \infty
\label{div2}
\end{align}

\noi
for any $t > 0$
and any non-zero test function $\psi \in \D(\T^2) = C^\infty(\T^2)$, 
where 
$\Xi_2\big((f_0, f_1)\big)$ denotes the second Picard iterate defined by 
 \begin{align}
 \Xi_2 \big((f_0, f_1)\big)(t)=  
 \int_0^t 
S(t- t') \NN\big(\dt S(t') f_0 + S(t') f_1 \big) dt'
\label{Pi3}
 \end{align}

\noi
with $S(t) = \frac{\sin(t\jb{\nb})}{\jb{\nb}}$.
See also \cite[Proposition 1.4]{Deya2} for a related result.

In a recent work~\cite{LLLOT}, 
 the authors with G.\,Li,  J.\,Li, and S.\,Liu
studied the quadratic NLW~\eqref{NLW1}
beyond the variance blowup
in the spirit of \cite{Hairer, LLOT}.
Fix $\be \ge \frac 12$.
Given $N \in \N$,
we introduce the multiplicative renormalization constant $\dl_{\be, N}$
by setting
\begin{align}
\dl_{\be, N}=
\begin{cases}
(\log \jb{N})^{-\frac 14} & \text{for } \be = \frac 12, \\
\jb{N}^{\frac 12 -\be}
  &  \text{for }   \be>\frac 12.
\end{cases}
\notag 
\end{align}

\noi
Note that 
$\dl_{\be, N} \to 0$ as $N \to \infty$.
Then,   consider the  quadratic NLW 
with the renormalized (frequency-truncated) random initial data:
\begin{align}
\begin{cases}
\dt^2 u_N  + (1 -  \Dl)  u_N = \NN(u_N)\\
(u_N, \dt u_N) |_{t = 0}
  = (\dl_{\be, N} \P_N u_0, \dl_{\be, N} \P_N u_1).
\end{cases}
\label{NLW2}
\end{align}

Our goal is to study the limiting behavior of $u_N$ as $N \to \infty$.
As in the BBM case, 
we first need to understand the convergence property of 
the second Picard iterate
$ \Xi_2(\dl_{\be, N} \P_N \vec u_0)$
for~\eqref{NLW2}, where 
$\Xi_2$ is as in \eqref{Pi3}
and 
$\vec u_0 = (u_0, u_1)$.
For $\frac 14 < \be < \frac 12$, it was
observed in~\cite[Proposition~1.4]{OO} that the second Picard
iterate $\Xi_2(\vec u_0 (\be))$ (without the renormalization constant $\dl_{\be, N}$) enjoys the bilinear dispersive smoothing
of $\frac 14$ as compared to the regularity one would expect from a naive parabolic thinking.
It turns out
 this bilinear dispersive smoothing is inherited 
 for $ \Xi_2(\dl_{\be, N} \P_N \vec u_0)$
 in the limit
as $N \to \infty$
even beyond $\be = \frac 12$.
Define the space-time noise
 $\eta$ on $\R_+ \times \T^2$ 
by the Fourier series: 
\begin{align*}
\eta(t) = \sum_{n\in \Z^2} \eta_n(t) e_n, 
\end{align*}

\noi
where  $\{\eta_n\}_{n \in \Z^2}$
is a family of  independent complex-valued
centered stationary Gaussian processes,
conditioned that
$\eta_{-n} = \cj{\eta_n}$, $n \in \Z^2$, 
with the  covariance given by 
\begin{align}
\E \Big[\eta_n(t_1) \cj{\eta_m(t_2)}\Big]
= \ind_{n = m} \cdot J_0\big((t_1-t_2)|n|\big).
\label{cov2}
\end{align}

\noi 
Here,  $J_0$ is the Bessel function of order zero.
It was shown in \cite{LLLOT}
that, as $N \to \infty$,  
the second Picard iterate
$ \Xi_2(\dl_{\be, N} \P_N \vec u_0)$
converges in law to 
\begin{align}
Z(t) = \int_0^t S(t - t') \eta (t') dt'
\label{cov3}
\end{align}

\noi
in $C([0, T]; W^{\frac 14 - \eps, \infty}(\T^2))$
for any 
$T > 0$ and $\eps > 0$.
Recall from 
\cite[Appendix B.8]{Gra}
that for $t_1 \ne t_2$, 
we have 
$\big|J_0\big((t_1-t_2)|n|\big)\big| \les_{t_1, t_2} |n|^{-\frac 12}$
as $|n| \to \infty$.
Hence, 
although the limit $Z$ in~\eqref{cov3}
is linear (which can not benefit from any multilinear dispersive smoothing), 
the covariance structure~\eqref{cov2}
given by the Bessel function $J_0$ encodes
the $\frac 14$-smoothing
(coming from the bilinear dispersive smoothing for each finite $N \in \N$).
As a result, 
the spatial regularity of $\eta(t)$ is $\frac 14$ higher
than that of a spatial white noise on $\T^2$.


We then 
proceed as in the BBM case
with 
the Skorokhod representation theorem \cite[Chapter~31]{Bass}
and work with the second order expansion for \eqref{NLW2}
on a new probability space.
We point out that 
the uniform (in $N$) solution theory for  \eqref{NLW2}
requires 
 multilinear dispersive 
smoothing and 
 the paracontrolled calculus as in \cite{GKO2}.
In \cite{LLLOT}, 
the authors  
proved
that, as $N \to \infty$,  
the (smooth) solution $u_N$ to \eqref{NLW2} converges 
almost surely (on a new probability space)
to the solution $u$ to 
the following stochastic NLW:
\begin{align} 
\begin{cases}
\dt^2  u + (1 - \Dl)  u =  u^2 + c_\be  \eta\\ 
( u, \dt  u)|_{t = 0} = (0, 0), 
\end{cases} 
\label{NLW3}
\end{align}

\noi
in $C([0, T]; H^{- \frac 12 - \eps}(\T^2))$
for any 
$\eps > 0$, 
where $c_\be > 0$ is some (deterministic) constant
 and $T = T_\o$ is an almost surely positive 
local existence time.
We note that the solution $u$
to the limiting equation \eqref{NLW3}
belongs to 
 $C([0, T]; H^{\frac 14 - \eps}(\T^2))$
and thus no renormalization is required
on the quadratic nonlinearity $u^2$ in \eqref{NLW3}.

We note that, due to the lack of global well-posedness
of the limiting equation~\eqref{NLW3}, 
the almost sure convergence 
of $u_N$ to $u$ 
on the new probability space
(after applying the Skorokhod representation theorem)
takes place on the random time interval $[0, T]$.
This is in contrast to the BBM case, 
where the convergence holds true on any finite interval $[0, T]$ 
with deterministic $T>0$.
Thus,  we need to proceed with care
in the current NLW case.

Define the real-valued random variables $A_N$ and $A$
by setting
\[A_N (\o) = \| u_N^\o \|_{C_{T_\o} H^{-\frac 12- \eps}_x} 
\qquad \text{and}\qquad 
A (\o)= \| u^\o\|_{C_{T_\o} H^{-\frac 12- \eps}_x}. \]

\noi
Note that the time interval $[0, T_\o]$ is random here.
Then, $A_N$ converges almost surely to $A$
(on the new probability space).
With this notation, 
%
it follows from  
the portmanteau theorem 
with the lower semi-continuity of the indicator
function on $(\ld, \infty)$
and 
the non-triviality of the solution $u$ to the limiting equation~\eqref{NLW3} that 
\begin{align}
\begin{split}
&  \limsup_{N \to \infty} 
\PP\Big(\| u_N\|_{C_T H^{-\frac 12- \eps}_x} >\ld\Big)\\
& \quad  \ge 
\liminf_{N \to \infty} 
\PP\Big(\| u_N\|_{C_T H^{-\frac 12- \eps}_x}>\ld\Big)
= \liminf_{N \to \infty} 
\E\big[\ind_{\{A_N >\ld\}}\big]\\
& \quad \ge  
\E\big[\ind_{\{A>\ld\}}\big]
 = \PP\Big(\| u\|_{C_T H^{-\frac 12- \eps}_x}> \ld\Big)
> 0 
\end{split}
\label{NLW4}
\end{align}

\noi
for some small $\ld > 0$.
Hence, from Proposition \ref{PROP:ill1}, 
we conclude the following probabilistic ill-posedness
 result
for the quadratic NLW \eqref{NLW1} on~$\T^2$.

\begin{theorem}
\label{THM:2}

Let $\be\geq \frac12$.
Then, 
the probabilistic Cauchy problem 
for the quadratic NLW~\eqref{NLW1}  on $\T^2$
with the Gaussian random initial data
$(u_0, u_1) = (u_0(\be), u_1(\be))$ in \eqref{gauss2}
is  ill-posed
in the sense of Definition~\ref{DEF:1}.

\end{theorem}

\begin{remark}

\rm

In \cite[Remark 4.3]{OO}, 
it was pointed out that 
variance blowup \eqref{div2}
for the quadratic NLW on the $d$-dimensional torus $\T^d$, $d \ge 3$, 
occurs 
for $\be \ge 1 - \frac d4$.
It would be of interest to investigate if 
 the ``beyond variance blowup'' result from \cite{LLLOT}
 and Theorem \ref{THM:2} extend to 
 higher dimensions $d \ge 3$.

\end{remark}

\begin{remark}\label{REM:NLS}\rm
In \cite{Liu}, 
R.\,Liu studied the well-posedness
issue of the following quadratic NLS on $\T^2$:\footnote{Here, 
we ignore the issue of renormalization.}
\begin{align}
i \dt u + \Dl u = |u|^2
\label{NLS3}
\end{align}

\noi
with the Gaussian random initial data $u_0(\al)$ in \eqref{gauss1}.
He proved that \eqref{NLS3}
is almost surely locally well-posed for $\al > \frac 12$, 
while variance blowup occurs for $\al \le \frac 14$.
It would be of interest to investigate 
what happens in the range $\frac 14 < \al \le \frac 12$
and to prove probabilistic ill-posedness as in Theorems \ref{THM:1}
and \ref{THM:2} 
for those values of $\al$, where variance blowup holds. 
See also Section \ref{SEC:3}.

\end{remark}

\subsection{On stochastic  PDEs}
\label{SUBSEC:SPDE}

In this subsection, 
we consider the following
stochastic dispersive PDE
with an additive forcing, posed on $\T^d$:\footnote{Once again, 
we ignore the issue of renormalization in the following discussion.}
\begin{align}
\begin{cases}
\dt u + Lu = \NN(u) +  \Phi \xi\\
u|_{t = 0} = u_0, 
\end{cases}
\label{SDE1}
\end{align}

\noi
where 
$L$ denotes a dispersive linear operator, 
$\NN(u)$ denotes 
a nonlinearity, 
and $\Phi$ is a (possibly unbounded) linear operator on $L^2(\T^d)$.
Here, 
$\xi$ denotes a space-time white noise 
on $\R_+\times \T^d$.
Stochastic dispersive PDEs
with additive noises have attracted much attention
and have been studied extensively;
see, for example, 
\cite{DD1, DDT1, Oh1, GKO, FOW, OO, ORW, GKOT, ORTz, GKO2, OQS, OOT1, BLL24, OOT2, 
LTzW, GOSW}.

We say that $u$ is a solution to \eqref{SDE1}
with initial data $u|_{t= 0} = u_0$
if it satisfies the following Duhamel formulation 
(= mild formulation):
\begin{align}
u(t) = e^{-tL } u_0 + \int_0^t e^{-(t-t')L} \NN(u)(t') dt'
+ \Psi(t), 
\notag
\end{align}

\noi
where $\Psi$ is the so-called stochastic convolution, 
 accounting for the effect of the
stochastic forcing in \eqref{SDE1}, 
and is formally  given by 
\begin{align}
\Psi(t) = \int_0^t e^{-(t-t') L} \Phi \xi (dt').
\label{SDE3}
\end{align}

\noi

Fix $s \in \R$.
As in Subsection~\ref{SUBSEC:1.2}, 
 we assume that 
\eqref{disp1}  (namely, \eqref{SDE1} without the stochastic forcing) is 
ill-posed
in $H^s(\T^d)$, while it is 
 locally well-posed
in $H^\s(\T^d)$ for some $\s > s$.
Our aim is to provide a brief discussion 
on the effect of rough noises and thus we assume that 
the (deterministic) initial data $u_0$ is sufficiently smooth
(for example, 
$u_0 \in H^\s(\T^d)$).
We also assume that the noise $\Phi \xi$ in \eqref{SDE1}
is spatially homogeneous.
Namely, $\Phi$ is a Fourier multiplier operator:
\[ \Phi(e_n) = \phi_n e_n, \quad n \in \Z^d.\]

\noi
For the sake of discussion, 
we further assume that $\Phi$ is Hilbert-Schmidt from $L^2(\T^d)$
to $H^s(\T^d)$
such that a standard argument yields that $\Psi \in C(\R_+; H^s(\T^d))$, almost surely.

Despite an additional difficulty coming from
the temporal roughness of the noise, 
the methodology
for proving local well-posedness
of \eqref{SDE1}
is analogous to that for proving 
almost sure local well-posedness
of \eqref{disp1}
with random initial data
as described in 
Subsection~\ref{SUBSEC:1.2}.
Namely, 

\smallskip
\begin{itemize}
\item[\bf Step 1:]
Starting with $\Xi_1 = \Xi_1(\Phi \xi): = \Psi$, 
we use stochastic analysis
and construct an enhanced data set 
$\Xi (\Phi\xi) = \big(\Xi_1(\Phi \xi), \dots, \Xi_K(\Phi \xi)\big)$, 
consisting of  random distributions
and random operators, 
and the space 
$ \cX_K$ for enhanced data sets
such that 
$\Xi \in \cX_K$, almost surely.

\medskip
\item[\bf Step 2:]
Given a (deterministic) enhanced data set $\Xi = (\Xi_1, \dots, \Xi_K) \in \cX_K$, 
write $u = \Xi_1 + \dots + \Xi_k + v$ 
for some $k \in \{1, \dots, K\}$ 
as in \eqref{exp1a}
and prove local well-posedness of
the equation 
for  the remainder term~$v$:
\begin{align}
\begin{cases}
\dt v + Lv = \NN(v + \Xi_1 + \dots + \Xi_k ) 
- (\dt + L) (\Xi_2 + \dots + \Xi_k ) \\
v|_{t = 0} = u_0
\end{cases}
\label{SDE3a}
\end{align}

\noi
by building  a continuous map:
\begin{align}
\Xi \in \cX_K
\longmapsto v \in C([0, T]; H^{\s}(\T^d))
\label{SDE3b}
\end{align}
 via deterministic analysis.

\end{itemize}

\smallskip

As in the random data case, 
by ``local well-posedness''\footnote{As is customary in the field, 
the term ``almost sure local well-posedness''
is used only for the random data case.
For a stochastic PDE, we simply use 
the term ``local well-posedness''.}
 of the stochastic dispersive PDE~\eqref{SDE1}, 
we usually mean that the following statements hold:

\smallskip

\begin{itemize}
\item[\bf (i)] \underline{\bf existence \& uniqueness.}  
There exist an almost surely positive stopping time $T$
and a process 
 $u$, almost surely belonging to  $C([0, T]; H^s(\T^d))$.
 Moreover, $u$ is a unique solution to \eqref{SDE1} on $[0, T]$, 
 almost surely.
As in the random data case,  both existence and uniqueness
of a solution are interpreted in a delicate sense, 
depending on a model under consideration.

\smallskip
\item[\bf (ii)]
\underline{\bf  stability under the frequency truncation.}
Given $N \in \N$, 
let 
$ u_N$
be the  solution  to 
\begin{align}
\begin{cases}
\dt u_N + Lu_N = \NN(u_N) +  \Phi \P_N \xi\\
u_N |_{t = 0} = u_0.
\end{cases}
\notag
\end{align}

\noi
Then, as $N \to \infty$, $u_N$ converges to 
$u$ 
 in $C([0, T];H^s(\T^d))$, almost surely, 
 where  $u$ 
 denotes the solution
to \eqref{SDE1}, 
constructed in Part (i).
Here,  $T$ denotes the 
almost surely positive  local existence time
 from Part (i) (possibly multiplied by a deterministic constant factor).\footnote{A similar
 comment on $T$ applies to (iii') below which we omit in the following.}

\end{itemize}

\smallskip

\begin{remark}\label{REM:sto2}\rm
As in the random data case (see Remark \ref{REM:sto1}), 
the construction of the enhanced data set $\Xi(\Phi\xi)$ in Step 1
is carried out by first 
 constructing
the enhanced data set  $\Xi(\Phi \P_N \xi )$
for the frequency-truncated noise $\Phi \P_N \xi $
and then 
by establishing 
almost sure convergence
of $\Xi(\Phi \P_N \xi )$ to a limit, denoted by $\Xi(\Phi \xi)$, 
in $\cX_K$.
Then, the stability under the frequency truncation 
follows from the continuity of the map~\eqref{SDE3b}
in Step 2.
\end{remark}

In the following, 
we focus on 
 the case of a monomial nonlinearity
 $  \NN(u) = u^p$ as in~\eqref{mono1}.
In this case,  a
modification of Step 1 typically 
 yields that, for each $k = 1, \dots, K$, 
there exists 
$\ld_k \in \N$ such that 
\begin{align}
\Xi_k(\dl \Phi\xi ) = \dl^{\ld_k} \Xi_k (\Phi\xi ) 
\notag
\end{align}

\noi
for any $\dl > 0$.
Namely, 
the enhanced data set 
 $\Xi(\dl  \Phi \xi )$, 
associated with  the equation:
\begin{align}
\dt u^\dl + Lu^\dl = \NN(u^\dl) + \dl  \Phi \xi, 
\label{SDE5}
\end{align}

\noi
converges almost surely to~$0$
in $\cX_K$
as $\dl \to 0$:
\begin{align}
\Xi(\dl\Phi\xi)= \big(\Xi_1(\dl \Phi\xi), \dots, \Xi_K(\dl \Phi\xi)\big)\,  \too \, \Xi (0) = (0, \dots, 0).
\label{SDE6}
\end{align}

\noi
This  leads to  the following stability in the small noise limit.

\smallskip

\begin{itemize}
\item[\bf (iii)] \underline{\bf stability 
in the small noise limit.}
Given $\dl > 0$, 
let 
$ u^\dl$
be the   solution  to~\eqref{SDE5} with $u^\dl|_{t = 0} = u_0$,
given by the expansion (see \eqref{exp1a}):
\begin{align}
u^\dl  = \Xi_1 (\dl \Phi\xi) + \dots + \Xi_k(\dl \Phi\xi) + v^\dl, 
\notag
\end{align}

\noi
where $v^\dl$ is a solution 
to \eqref{SDE3a} with $\Xi = \Xi(\dl \Phi \xi)$.
As in the random data case, 
we can choose the almost surely positive local existence  time $T = T_\o$
to be independent of $0  < \dl \le 1$.
Then, it follows from  \eqref{SDE6} with Step 2 above
that, as $\dl \to 0$, 
\begin{align}
 u^\dl \too  u
\label{SDE7}
\end{align}

\noi
 in $C([0, T];H^s(\T^d))$, almost surely, 
 where $u$ is the solution to the deterministic equation~\eqref{disp1}
 with the same initial data $u|_{t = 0} = u_0$.

\end{itemize}

\smallskip

The stability in the small noise limit
is of importance, 
guaranteeing
that a weak stochastic forcing leads
to a mild perturbation rather than 
generating erratic behavior, 
and is a basis of the study on the large deviation principle
as $\dl \to 0$ (see \cite[Theorem 5.7]{BGOR}; see also \cite{CX})
and on metastability (see \cite{BL}).

Arguing as in \eqref{disp7a} and \eqref{disp9}
with Remark \ref{REM:sto2}, 
we arrive at the following version
of the stability in the small noise limit.

\smallskip

\begin{itemize}
\item[\bf (iii')]
\underline{\bf stability in the small noise limit.}
Given 
 a sequence 
  $\{\dl_N\}_{N \in \N}$ 
   of positive numbers tending to $0$
 as $N \to \infty$, 
 let 
$ u_N$
be  the   solution  to 
\begin{align}
\begin{cases}
\dt u_N + Lu_N = \NN(u_N) + \dl_N  \Phi \P_N \xi\\
u_N |_{t = 0} = u_0.
\end{cases}
\label{SDE8}
\end{align}
 
\noi
Then,  as $N \to \infty$, we have 
\begin{align}
 u_N \too  u
\label{SDE9}
\end{align}

\noi
 in $C([0, T];H^s(\T^d))$, almost surely, 
 where $u$ is the solution to the deterministic equation~\eqref{disp1}
 with initial data $u|_{t = 0} = u_0$.

\end{itemize}

\smallskip

The stability conditions (iii) and (iii')
are {\it not} for a single equation \eqref{SDE1}
but for the family of the equations \eqref{SDE5}
or \eqref{SDE8}.
As such, 
the failure of the condition (iii) or~(iii')
does not lead to ``ill-posedness'' of the original stochastic dispersive PDE~\eqref{SDE1}.
Nonetheless, by noting that the 
convergence \eqref{SDE7} and \eqref{SDE9}
to the deterministic solution $u$
in the small noise limit
is an SPDE analogue of the law of large numbers, 
we conclude that
the failure of the condition (iii) or (iii')
 implies that, 
even in the small noise limit, 
the effect of the stochastic forcing does not become negligible, 
preventing 
a ``law of large numbers''-type result
(and also large deviation and metastability results)
for the family of the equations \eqref{SDE5}.

\medskip

The ``beyond variance blowup'' results
discussed in Subsections \ref{SUBSEC:2.1} 
and  \ref{SUBSEC:2.2} 
have direct analogues
in the stochastically forced setting, 
violating the stability condition (iii') above, 
which we briefly discuss below.

Given $\be \in \R$, 
consider the following stochastic BBM:
\begin{align}
\dt u-\partial_{txx}u+\partial_{x}u
+\partial_{x}(u^{2})  = \jb{\dx}^{\be} \dx \xi,   
\label{BBM6}
\end{align}

\noi
where $\xi$ denotes the space-time white noise on $\R_+\times \T$.
It follows from a slight modification of the argument in \cite{Forlano}
that 
\eqref{BBM6} is locally well-posed for $\be < \frac 34$.
Moreover, variance blowup occurs
for $\be \ge \frac 34$.

Let    $\be \ge \frac 34$.
Given a mean-zero function $u_0 \in H^1(\T)$, 
 consider the following renormalized stochastic BBM:
\begin{align}
\begin{cases}
\dt u_N-\partial_{txx}u_N+\partial_{x}u_N
+\partial_{x}(u_N^{2})  = \dl_{1-\be, N}  \jb{\dx}^\be  \dx \P_N \xi\\
u_N |_{t = 0} = u_0, 
\end{cases}
\label{BBM7}
\end{align}

\noi
where 
$\dl_{1-\be, N}$ is as in \eqref{CN}.
Then, 
it was shown in 
\cite{LLOT} that,  
as $N \to \infty$, 
the solution $u_N$
to~\eqref{BBM7}
with the renormalized noise
converges in law
to $u$
in $C(\R_+; H^{- \frac 14 - \eps}(\T))$
for any $\eps > 0$, 
where $u$ is the solution to 
the following stochastic BBM:
\begin{align*}
\begin{cases}
\dt u-\partial_{txx}u+\partial_{x}u
+\partial_{x}(u^{2})  = \dx \upze  \\
u|_{t = 0} = u_0.
\end{cases}
\end{align*}

\noi
Here, 
the space-time noise
 $\upze$ is given by 
\begin{align*}
\upze (t)  = \sum_{n\in \Z^*}
\int_0^t \sqrt {2t'} d B_n(t') e_n
\end{align*}

\noi
for each $t \in \R_+$,  
where 
 $\{ B_n \}_{n \in \Z}$ is a family of mutually independent complex-valued
Brownian motions with $\text{Var}(B_n(t)) = t$ conditioned  that $B_{-n} = \cj{B_n}$, $n \in \Z$. 
In particular, for each $t \in \R_+$, 
 $\upze(t)$ is a scalar multiple of a spatial
white noise 
 with variance  $t^2$
(with spatial mean~$0$). 

An analogous ``beyond variance blowup'' result
holds for the following stochastic NLW on $\T^2$
with $\be \ge \frac 12$:
\begin{align}
\dt^2 u + (1- \Dl) u = \, :\!u^2\!:  + \, \jb{\nb}^\be \xi, 
\label{SNLW1}
\end{align}

\noi
where $\xi$ denotes the space-time white noise on $\R_+\times \T^2$.
See \cite{LLLOT} for details.
These results show that the stability condition (iii')
fails for

\smallskip

\begin{itemize}
\item  the stochastic BBM \eqref{BBM6} with $\Phi = \jb{\dx}^{\be} \dx$ and $\be\ge  \frac 34$, 

\medskip

\item  the stochastic NLW \eqref{SNLW1} on $\T^2$ with $\Phi = \jb{\nb}^{\be} $ and $\be \ge \frac 12$.

\end{itemize}

\begin{remark}\rm 
In the examples above, 
the linear operator $\Phi$ in \eqref{SDE1} took a specific form:
 $\Phi = \jb{\dx}^{\be} \dx$ for the stochastic BBM \eqref{BBM6}
 and  $\Phi = \jb{\nb}^{\be} $ for the stochastic NLW \eqref{SNLW1}.
Suppose that, given general $\Phi \in \HS(L^2(\T^d); H^s(\T^d))$
for some $s \in \R$,   the stochastic convolution $\Psi(t)$ in \eqref{SDE3}
is merely a spatial distribution of negative regularity for each fixed $t > 0$.
In this case, in order to construct a Wick-renormalized power
$:\!\Psi^k\!:$\,, 
we need to additionally assume that 
$\Phi$ is $\g$-radonifying from $L^2(\T^d)$ to a certain Fourier-Lebesgue space.
See 
\cite{OOPTz} for a further discussion.
This comment also applies to the parabolic case;
see \cite[Remark 1.5]{OOPTz}.

\end{remark}

We conclude this section by extending the discussion above
to stochastic parabolic PDEs.
In the study of singular stochastic parabolic PDEs, 
 a  variance blowup phenomenon
has been observed for the following fractional KPZ equation on $\T$:\footnote{For simplicity
of the presentation, we ignore the renormalization ``$(\dx h)^2 - \infty$''
on the quadratic nonlinearity $(\dx h)^2$ in \eqref{KPZ1}.}
\begin{align}
\dt h = \dx^2 h  + (\dx h)^2 + |\dx|^{\be}\xi
\label{KPZ1}
\end{align}

\noi
for  $\be \ge \frac 14$, where $\xi$ denotes a space-time white noise
on $\R_+\times \T$;
see \cite[Subsection 4.9]{Hoshino}.
In~\cite{Hairer}, 
Hairer 
considered the following renormalized fractional 
KPZ (with $\be = 1$):
\begin{align}
\dt h_N = \dx^2 h_N + \big((\dx h_N)^2 - C_N \big) + N^{-\frac 34}\P_N  \dx \xi, 
\label{KPZ2}
\end{align}

\noi
where
the frequency-truncated
noise $\P_N  \dx \xi$
is endowed with 
a vanishing multiplicative renormalization constant
$ N^{-\frac 34}$.
Then, 
under  an appropriate assumption on deterministic initial data, 
he showed that, as $N \to \infty$,  the solution $h_N$  to \eqref{KPZ2}
converges  in law
to a solution to the standard KPZ forced by a space-time white noise
(namely, \eqref{KPZ1} with $\be = 0$ but with a suitable variance), 
thus violating (a suitable analogue of)
the stability in  the small noise limit (iii').

\begin{remark}\rm
As pointed out above, the stability in  the small noise limit (iii')
is 
a ``law of large numbers''-type statement.
On the other hand, 
as mentioned 
in Subsection \ref{SUBSEC:2.1}, 
the ``beyond variance blowup'' results
in \cite{Hairer, LLOT, LLLOT}
are
``central limit theorem'' type results.
Namely, they belong to completely different regimes.

\end{remark}

\section{Variance blowup 
as mild probabilistic ill-posedness}
\label{SEC:3}

In this section, we provide a brief discussion
on how to view variance blowup results
as the failure of smoothness at the origin
in the probabilistic sense.

In \cite[Section 6]{BO97}, 
Bourgain considered the (deterministic) well-\,/\,ill-posedness
issue for the Korteweg-de Vries equation (KdV) in the low regularity setting:
\begin{align}
\dt u + \dx^3 u + \dx (u^2) = 0
\label{KDV1}
\end{align}

\noi
by studying the smoothness of the solution map
at the origin.
In particular, 
he showed that the solution map $\G:
u_0 \in H^s(\M) \mapsto u \in C([0, T]; H^s(\M))$, 
a priori defined on smooth initial data, 
fails to be $C^3$ at the origin
for 
\begin{align}
s < - \frac 12 \ \text{ when } \M = \T
\quad \text{and}\quad 
s < - \frac 34 \ \text{ when }  \M = \R.
\label{KDV2}
\end{align}
See also~\cite{Tz99}.
This failure of smoothness
of the solution map
does not imply ill-posedness of KdV \eqref{KDV1}
in the regularity regime \eqref{KDV2}.
However, it 
implies 
{\it mild ill-posedness} of \eqref{KDV1}
in the regularity regime \eqref{KDV2}
in the sense 
that a contraction argument
can not be used to prove local well-posedness of \eqref{KDV1}
under \eqref{KDV2}, 
since a contraction argument would
provide analyticity\footnote{Hereafter, analyticity simply means
real analyticity.} of the solution map.
This result was strengthened to the failure 
of local uniform continuity of the solution map in \cite{CCT}
(with artificial lower bounds on $s$).
We note that 
the results \cite{BO97, Tz99, CCT}
complements
the local well-posedness result of \eqref{KDV1}
in 
\cite{KPV}
(modulo the endpoint).

\begin{remark}\rm
Strictly speaking, in the periodic setting, 
we need to restrict
our attention to the subclass $H^s_0(\T)$
of $H^s(\T)$, consisting of mean-zero functions, 
in the discussion above.
Without such a restriction, 
the solution map for KdV \eqref{KDV1} on $\T$
is {\it not} locally uniformly continuous for any $s \in \R$
(and thus is not $C^k$-smooth for any $k \in \N$).
See, for example,  
\cite[Remark 1.2 and Proposition 7.1]{Moli}
and 
\cite[Theorem~3.1.1 and Proposition~3.2.2]{Herr}.
We, however, ignore this issue in the following.

\end{remark}

Suppose that the solution map $\G:
u_0 \in H^s(\M) \mapsto u \in C([0, T]; H^s(\M))$  is 
$C^k$-smooth.
Then, given  smooth $u_0 \in H^\infty(\M)$
(guaranteeing
existence
of the solution $\G(\dl u_0)$  to \eqref{KDV1}
with initial data $\dl u_0$
on a time interval $[0, T]$, uniformly in $|\dl |\le 1$), 
it follows from the chain rule that 
\begin{align}
\bigg\|\frac {d^k}{d\dl^k} \G(\dl u_0)\Big|_{\dl = 0} \bigg\|_{C_{T} H^s_x}
\les \| u_0\|_{H^s}^k.
\label{mild0}
\end{align}

\noi
Bourgain's argument in \cite{BO97} is based
on 
considering
the (smooth) solution $u_{N}^\dl = \G(\dl u_{0, N})$
to~\eqref{KDV1} with initial data
$u_{N}^\dl|_{t = 0} = 
\dl u_{0, N}$
for  a suitable sequence $\{u_{0, N}\}_{N  \in \N}$
of smooth functions
$u_{0, N}$
(with finite Fourier support)
such that $\|u_{0,N}\|_{H^s} = O(1)$, 
and showing that 
$\dd_\dl^3 u_{N}^\dl|_{\dl = 0}$
is unbounded in $C([0, T] ; H^s(\M))$,
violating \eqref{mild0} with $k = 3$, 
which 
 implied the failure
of $C^3$-smoothness
of the solution map at the origin.

\medskip

Let us now turn to the random data setting
and consider the model dispersive PDE~\eqref{disp1}
on $\T^d$ 
with a monomial nonlinearity \eqref{mono1} which is analytic
in the unknown $u$ (and its conjugate in the complex-valued case).

Suppose that, given $\al \in \R$, 
\eqref{disp1} is almost surely locally well-posed
with the Gaussian random initial data $u_0 = u_0(\al)$ in \eqref{gauss1}
in the sense described in 
Subsection \ref{SUBSEC:1.2}.
Let $\cX_K$ be the space of enhanced data sets 
$\Xi = (\Xi_1, \dots, \Xi_K)$ determined in Step 1
such that $\Xi(u_0^\o) = \big(\Xi_1(u_0^\o), \dots, \Xi_K(u_0^\o)\big) \in \cX_K$
for each $\o \in \Si$.
Here,  $\Si\subset \O$
is the set of full probability, 
where  the almost sure local well-posedness of \eqref{disp1} holds.
Furthermore, 
we  assume that, 
 local well-posedness
of the equation \eqref{disp2} for the remainder term $v$ under the expansion \eqref{exp1a}
follows
from a contraction argument,  
which in particular implies that 
the map \eqref{disp4}
is analytic in $\Xi = (\Xi_1, \dots, \Xi_K)$.

Now, given 
$\dl \in \R$ with  $|\dl| \le 1$, 
let $u^\dl$ be the random solution to  \eqref{disp1}
with the scaled random initial data $u^\dl|_{t = 0} = \dl u_0$.
Then, under the assumption \eqref{disp5a}, 
we see that, for each $\o \in \Si$,  the map:
\[\dl \longmapsto 
\Xi(\dl u_0^\o)= \big(\Xi_1(\dl u_0^\o), \dots, \Xi_K(\dl u_0^\o)\big)\]

\noi
is analytic in $\dl$.
Then, combining this with the analyticity of the map \eqref{disp4}, 
we obtain the following 
smoothness  at the origin.

\smallskip

\begin{itemize}
\item[\bf (iv)] \underline{\bf smoothness  at the origin.}
With the assumption and notation above, 
we have that, 
 for each $\o \in \Si$,  the map:
\[ \dl \longmapsto 
u^{\dl, \o} = 
\Xi_1(\dl u_0^\o)
+ \dots + \Xi_k(\dl u_0^\o) + v^{\dl, \o} \in C([0, T_\o]; H^s(\T^d))\]

\noi
is analytic in $|\dl| \le 1$.

\end{itemize}

\smallskip

Given $N \in \N$, 
let $ u_N^\dl = \G(\dl \P_{N} u_0)$ be
 the  (smooth) random solution  to \eqref{disp1} with the scaled frequency-truncated
 random initial data
$ u_N^\dl|_{t = 0} = \dl \P_{N} u_0$.
An analogous discussion 
together with Remark~\ref{REM:sto1}
leads to the following version 
of smoothness at the origin.

\smallskip

\begin{itemize}
\item[\bf (iv')] \underline{\bf smoothness  at the origin.}
For each $\o \in \Si$,  the map:
\[ \dl \longmapsto 
u^{\dl, \o}_N 
 \in C([0, T_\o]; H^s(\T^d))\]

\noi
is analytic in $|\dl| \le 1$, 
uniformly in $N \in \N$
in the sense that,  
 for each $k \in \N$ and $\o \in \Si$, we have 
\begin{align}
\sup_{N \in \N}
\bigg\|\frac {d^k}{d\dl^k} u_N^{\dl, \o} \Big|_{\dl = 0} \bigg\|_{C_{T_\o} H^s_x}
= 
\sup_{N \in \N}
\bigg\|\frac {d^k}{d\dl^k} \G(\dl \P_{N} u_0^\o) \Big|_{\dl = 0} \bigg\|_{C_{T_\o} H^s_x}
\le C_{k, \o} < \infty.
\label{mild1}
\end{align}

\end{itemize}

\smallskip

The condition \eqref{mild1} for fixed $k \in \N$
can be viewed as a probabilistic analogue of \eqref{mild0}, 
leading to the following definition
of mild probabilistic ill-posedness.

\begin{definition}\label{DEF:2}\rm
Given $\al \in \R$,
we say that 
the probabilistic Cauchy problem \eqref{disp1} 
with the Gaussian random initial data
$u_0 = u_0(\al)$ in~\eqref{gauss1}
is  mildly ill-posed
for this particular value of $\al$
if the condition \eqref{mild1} fails for some integer $k \ge 2$.

\end{definition}

When $k = 1$, we have
\[ \sup_{N \in \N} \big\| \dd_\dl u_N^{\dl, \o}|_{\dl = 0} \big\|_{C_{T_\o} H^s_x}
= \sup_{N \in \N}  \|e^{-tL} \P_N u_0^\o\|_{C_{T_\o} H^s_x}
\le  \|u_0^\o\|_{ H^s}, 
\]

\noi
trivially satisfying \eqref{mild1}.
This is the reason that 
the value $k = 1$ was excluded in Definition~\ref{DEF:2}.

In the following, we assume that 
\eqref{disp1} is locally well-posed in $H^\s(\T^d)$ via
a contraction argument.
Then, given any $f \in H^\s(\T^d)$, 
the solution $u$ to \eqref{disp1}
with initial data $u|_{t = 0} = f$
is given as the limit of Picard iterations
(or by a power series expansion).
Namely, there exists a family $\{A_j\}_{j \in \N}$
of $m_j$-linear operators (for some increasing sequence $\{m_j\}_{j \in \N} \subset \N$
with $m_1 = 1$), 
bounded from $\big(H^\s(\T^d)\big)^{\otimes m_j}$
into $C([0, T]; H^\s(\T^d))$ (where $T = T(f) > 0$ denotes
the local existence time), 
such that 
\begin{align}
u = u(f) = \sum_{j = 1}^\infty A_j (f, \dots, f) .
\notag
\end{align}

\noi
In particular, 
for the solution  $ u_N^\dl$ 
 to \eqref{disp1} with the scaled frequency-truncated
 random initial data
$ u_N^\dl |_{t = 0} = \dl \P_{N} u_0$, we have
\begin{align}
\begin{split}
u^\dl_N  
& = \sum_{j = 1}^\infty A_j (\dl \P_{N} u_0, \dots, \dl \P_{N} u_0) \\
& = \sum_{j = 1}^\infty \dl^{m_j} A_j ( \P_{N} u_0, \dots,  \P_{N} u_0), 
\end{split}
\notag
\end{align}

\noi
where the second equality follows
from the multilinearity of $A_j$.
Hence, we have 
\begin{align}
\frac {d^k}{d\dl^k} u_N^\dl \bigg|_{\dl = 0}
= \ind_{k = m_j}
\cdot m_j !
\cdot 
A_j ( \P_{N} u_0, \dots,  \P_{N} u_0)
\notag
\end{align}

\noi
for $k \in \N$.

Note that when $j = 2$, 
the term $A_2 ( \P_{N} u_0,  \P_{N} u_0)$ corresponds to the second Picard iterate:
\begin{align}
A_2 ( \P_{N} u_0,  \P_{N} u_0)(t) 
= \I \big( \NN(\P_N u_0) \big)(t), 
\notag
\end{align}

\noi
where $\I$ is the Duhamel integral operator defined in 
\eqref{disp4b}.
Then, the variance blowup 
on the (random) second Picard iterate such as \eqref{div1}
and \eqref{div2}
implies that \eqref{mild1} fails for $k = 2$, 
implying
mild probabilistic ill-posedness
in the sense of Definition \ref{DEF:2}.

%
%

Let us conclude this paper by discussing the case of 
the 
renormalized quadratic NLS on~$\T^2$:
\begin{align}
i \dt u + \Dl u = \NN(u),  
\label{NLS4}
\end{align}

\noi
where $\NN(u)$ denotes the renormalized nonlinearity given  by
\[ \NN(u) = \P_{\ne 0}(|u|^2) = |u|^2 - \int_{\T^2} |u|^2 dx.\]

\noi
Here, $\P_{\ne 0}$ denotes the projection onto non-zero (spatial) frequencies.
In \cite{LO}, 
sharp local well-posedness of \eqref{NLS4} in $L^2(\T^2)$
was established via a contraction argument.

Consider \eqref{NLS4} 
with the Gaussian random initial data $u_0 = u_0(\al)$
in \eqref{gauss1}.
As mentioned in Remark \ref{REM:NLS}, 
R.\,Liu proved
almost sure local well-posedness of \eqref{NLS3}
for $\al > \frac 12$, 
while, for $\al \le \frac 14$, 
he proved the following variance blowup result:
\begin{align}
\lim_{N \to \infty} 
\E\Big[|\jb{  \Xi_2(\P_N u_0)(t), e_n}|^2\Big] = \infty
\label{div3}
\end{align}

\noi
for $t > 0$ and $n \in \Z^2\setminus\{0\}$, 
where 
$\Xi_2(f)$ denotes the second Picard iterate defined by 
 \begin{align*}
 \Xi_2 (f)(t)=  
 \int_0^t 
e^{i (t- t') \Dl} \NN(e^{it'\Dl } f ) dt'.
 \end{align*}

\noi
See \cite[Proposition 1.5]{Liu}.
In view of the discussion above, 
the variance blowup \eqref{div3}
implies the failure of the bound \eqref{mild1}
for $k = 2$.
Hence, we conclude that 
the 
renormalized quadratic NLS~\eqref{NLS4}
 on~$\T^2$ with the Gaussian random initial data 
 $u_0 = u_0(\al)$ in \eqref{gauss1}
is mildly probabilistically ill-posed for $\al \le \frac 14$.

\begin{ackno}\rm
The authors
would like to thank Shao Liu for his careful reading of the paper.
T.O.~was supported by the European Research Council (grant no.~864138 ``SingStochDispDyn'')
and also 
acknowledges support from  
the NSFC (grant no.~W2531005).
N.T. was partially supported by the ANR
project Smooth ANR-22-CE40-0017.

\end{ackno}

\end{document}